\documentclass{sl}

\usepackage[numbers]{natbib}
\usepackage{latexsym}
\usepackage{ae,aecompl} 
\usepackage{url}
\usepackage{hyperref}
\usepackage{amsmath}
\usepackage{amssymb}

\usepackage{slsec,stud_log,slfoot,slthm}

\makeatletter
\def\ps@SLppKAP{\let\@mkboth\@gobbletwo
\def\@oddhead{}\def\@oddfoot{\scriptsize\it\vbox{\baselineskip 9pt%
\noindent{\rmVIIhalf forthcoming in \itVIIhalf Studia Logica} %
{\rmVIIhalf (\therok).\newline
\raise 1pt\hbox{\tiny\copyright} \therok} %
{\itVIIhalf Kluwer Academic Publishers.\ %
Printed in the Netherlands.}\hfill}}%
\def\@evenhead{}\def\@evenfoot{}\def\chaptermark##1{}\def\sectionmark##1{}%
\def\subsectionmark##1{}}
\makeatother

\newif\ifeswirdpdf
\makeatletter
\@ifundefined{pdfoutput}%
  {
    \eswirdpdffalse
  }%
  {\ifcase\pdfoutput
      \eswirdpdffalse
   \else
    \eswirdpdftrue
   \fi}%
\makeatother

\ifeswirdpdf
\else
  \usepackage{aeguill}  
\fi

\newtheorem{theorem}{Theorem}
\newtheorem*{theorem*}{Theorem}
\newtheorem{proposition}[theorem]{Proposition}
\newtheorem{lemma}[theorem]{Lemma}

\newtheorem{corollary}[theorem]{Corollary}
\theoremstyle{definition}       
\newtheorem{definition}[theorem]{Definition}

\newtheorem{example}[theorem]{Example}

\let\meps\varepsilon
\let\epsilon\varepsilon 
\def\ar{\mathrm{ar}}
\def\deg{\mathrm{deg}}
\def\rk{\mathrm{rk}}

\def\PC{\ensuremath{\mathrm{PC}}}
\def\EC{\ensuremath{\mathrm{EC}}}
\def\AxEC{\mathrm{AxEC}}
\def\AxPC{\mathrm{AxPC}}
\newcommand*{\ECeps}{\EC_{\meps}}
\newcommand*{\PCeps}{\PC_{\meps}}
\def\Ax{\mathrm{Ax}}
\newcommand*{\tpkt}{\rlap{$\;$.}}
\newcommand*{\tkom}{\rlap{$\;$,}}
\newcommand*{\eps}[2]{\varepsilon_{#1} #2}
\newcommand*{\seq}{\mathrel{\Rightarrow}}
\newcommand*{\impl}{\mathrel{\rightarrow}}
\newcommand*{\cc}[1]{\mathrm{cc}(#1)}
\newcommand*{\card}[1]{\ensuremath{\lvert {#1} \rvert}}
\newcommand*{\size}[1]{\mathrm{sz}(#1)}
\DeclareMathOperator{\width}{\mathrm{wd}}
\DeclareMathOperator{\expan}{\mathrm{len}}
\newcommand*{\HC}{\mathrm{HC}}

\def\Hyp{\mathrm{Hyp}}

\begin{document}
\setcounter{page}{1} \thispagestyle{SLppKAP} \label{p}

\twoAuthorsTitle{G.\ts Moser}{R.\ts Zach}{The Epsilon Calculus and\\
  Herbrand Complexity}

\Abstract{Hilbert's $\meps$-calculus is based on an extension of the language
  of predicate logic by a term-forming operator $\meps_{x}$.  Two fundamental
  results about the $\meps$-calculus, the first and second epsilon theorem,
  play a r\^ole similar to that which the cut-elimination theorem plays in
  sequent calculus. In particular, Herbrand's Theorem is a consequence of the
  epsilon theorems.  The paper investigates the epsilon theorems and the
  complexity of the elimination procedure underlying their proof, as well as
  the length of Herbrand disjunctions of existential theorems obtained by this
  elimination procedure.}


\Section{Introduction}

Hilbert's $\meps$-calculus \cite{HB70,Leisenring:1969} is based on an
extension of the language of predicate logic by a term-forming operator
$\meps_{x}$.  This operator is governed by the \emph{critical axiom}
\begin{equation*}
 A(t) \impl A(\eps{x}{A(x)}) \tkom
\end{equation*}
where $t$ is an arbitrary term. Within the $\meps$-calculus, quantifiers
become definable by $\exists x A(x) \Leftrightarrow A(\eps{x}{A(x)})$ and
$\forall xA(x) \Leftrightarrow A(\eps{x}{\lnot A(x)})$.

The $\meps$-calculus was originally developed in the context of Hilbert's
program of consistency proofs. Early work in proof theory (before Gentzen)
concentrated on the $\meps$-calculus and the $\meps$-substitution method and
was carried out by Ackermann~\cite{Ack25,Ack40} (see
also~\cite{Moser:APAL:2005}), von Neumann~\cite{Neu27}, and Bernays
\cite{HB70} (see also \cite{Zach:2002,Zach:2004}).  The $\meps$-calculus is of
independent and lasting interest, however, and a study from a computational
and proof-theoretic point of view is particularly worthwhile.

The aim of this paper is to present the central notions of and basic results
about the $\meps$-calculus in a streamlined form and with attention to
questions of proof complexity.  This, we hope, will make the $\meps$-calculus
more easily accessible to a broader audience, and will make clearer the merits
and disadvantages of the $\meps$-calculus as a formalization of first-order
logic.

One simple example of the merits of the $\meps$-calculus is that by
encoding quantifiers on the term-level, formalizing (informal) proofs is
sometimes easier in the $\meps$-calculus, compared to formalizing proofs in,
e.g., sequent calculi. This is possible as the $\meps$-calculus allows more
condensed representation of proofs than standard sequent- or natural
deduction calculus.  However, it should be pointed out that the encoding of
quantifiers on the term-level can come at a significant cost, as the
transformation of quantified formulas may result in rather complicated
term-structures.

Hilbert's $\meps$-calculus is primarily a classical formalism, and we will
restrict our attention to classical first-order logic. (But see the work
of Bell~\cite{Bell:93a,Bell:93b}, DeVidi~\cite{DeVidi:95},
Fitting~\cite{Fitting:75}, Mostowski~\cite{Mostowski63} for a complementary
view.)  Our study is also motivated by the recent renewed interest in the
$\meps$-calculus and the $\meps$-substitution method in, e.g., the work of
Arai~\cite{Arai01b,Arai:2005}, Avigad~\cite{Avigad01}, and Mints et al.,
\cite{Mints:2003,Mints:99}.  The $\meps$-calculus also allows the
incorporation of choice construction into logic~\cite{BlassGurevich:2000:JSL}.

Our main focus will be the presentation and analysis of two
\emph{conservativity results} for the $\meps$-calculus: the first and second
$\meps$-theorems. Our proofs of these results are essentially those due to
Bernays~\cite{HB70}, but we are also concerned with bounds on the length of
proofs and the length of Herbrand disjunctions.  Let $T$ denote a finitely
axiomatized open theory with axioms $P_1$, \dots, $P_t$ containing no
quantifiers or $\meps$-terms, and let $\PCeps$ be a usual formulation of the
predicate calculus extended by the $\meps$-operator and its characteristic
axiom.  Then the first $\meps$-theorem states that any formula without
quantifiers or $\meps$-terms provable in $\PCeps$ from $T$ is already provable
from $T$ in the quantifier- and $\meps$-free fragment $\EC$ of $\PC$ (the
so-called elementary calculus of free variables).  The second $\meps$-theorem
says that any formula without $\meps$-terms provable in $\PCeps$ from $T$ is
also provable from $T$ in $\PC$, the pure predicate calculus.

We prove the first $\meps$-theorem in Section~5.  An important extension of
the first $\meps$-theorem is the so-called extended first $\meps$-theorem,
which yields one direction of Herbrand's theorem: if $A$ is an existential
formula provable in $\PCeps$, then there is a disjunction of instances of the
matrix of~$A$ provable in~$\EC$.  Call the minimal length of such a
disjunction the \emph{Herbrand complexity}~$\HC(A)$ of $A$. We obtain an upper
bound on $\HC(A)$ by analyzing the complexity of the procedure used to
eliminate critical formulas from the proof of~$A$; this bound is
hyperexponential in the number of quantifier axioms in the derivation of~$A$.
In Section~6 we show that there is also a hyperexponential lower bound
on~$\HC(A)$, and thus the procedure for generating the Herbrand disjunction
based on Bernays's proof of the extended first $\meps$-theorem is essentially
optimal.  In Section~7, we prove the second $\meps$-theorem.  Its proof also
contains a proof of the second direction of Herbrand's theorem, i.e., that
given a (valid) Herbrand disjunction of $A$, $A$ is provable.

In this paper we only consider the $\meps$-calculus without equality.  The
case of the $\meps$-theorems with equality is much more involved than the
proofs we consider here, and cannot be dealt with adequately in the space
available.  It raises important and interesting issues for the suitability of
formalisms based on the epsilon calculus for automated theorem proving.  The
complexity of the epsilon theorems in the presence of equality will be the
subject of future work.

\Section{The Epsilon Calculus: Syntax}\label{sec:syntax}

The syntax of the epsilon calculus is essentially that of a standard
first-order language.  We will, however, frequently pass back and forth
between different calculi formulated in slightly different languages.  Let us
first dispense with some pedantry: The language $L(\EC)$ of the elementary
calculus consists of: (1) free variables: $a$, $b$, $c$, \dots, (2) bound
variables: $x$, $y$, $z$, \dots, (3) constant and function symbols: $f$, $g$,
$h$, \dots with arities $\ar(f)$, \dots, (4) predicate symbols: $P$, $Q$, $R$,
\dots with arities $\ar(P)$, \dots, and the propositional connectives:
$\land$, $\lor$, $\impl$, $\lnot$.  The language $L(\PC)$ of the predicate
calculus is $L(\EC)$ plus the quantifiers $\forall$, $\exists$.  The language
$L(\ECeps)$ of the pure epsilon calculus is $L(\EC)$ plus the epsilon
operator~$\meps$.  The language $L(\PCeps)$ of the predicate calculus with
epsilon is $L(\PC)$ together with $\meps$.

We will distinguish between terms and semi-terms, and between formulas and
semi-formulas.  The definitions are:

\begin{definition} \label{d:termsformulas}
  \emph{(Semi)terms}, \emph{(semi)formulas}, and \emph{sub-(semi)terms} are
  defined as follows:
\begin{enumerate}
\item Any free variable $a$ is a (semi)term. Its only sub-(semi)term is $a$
  itself.  It has no immediate sub-(semi)terms.
\item Any bound variable $x$ is a semi-term. It has no sub-terms or immediate
  sub-(semi)terms. Its only sub-semiterm is $x$ itself.
\item If $f$ is a function symbol with $\ar(f) = 0$, then $f$ is a (semi)term.
Its only sub-(semi)term is $f$ itself.  It has no immediate sub-(semi)terms.
\item If $f$ is a function symbol with $\ar(f) = n > 0$, and $t_1$, \ldots,
  $t_n$ are (semi)terms, then $f(t_1, \ldots, t_n)$ is a (semi)term.  Its
  immediate sub-semiterms are $t_1$, \dots, $t_n$, and its immediate
  sub-terms are those among $t_1$, \dots, $t_n$ which are terms, plus the
  immediate subterms of those among $t_1$, \ldots, $t_n$ which are not terms.
  Its sub-semiterms are $f(t_1, \ldots, t_n)$ and the sub-semiterms of
  $t_1$, \dots, $t_n$; its subterms are those of its sub-semiterms which are
  terms.
\item If $P$ is a predicate symbol with $\ar(P) = n > 0$, and $t_1$,
\dots, $t_n$ are (semi)terms, then $P(t_1, \ldots, t_n)$ is an (atomic)
(semi)formula. Its
immediate sub-semiterms are $t_1$, \dots, $t_n$.  Its immediate
subterms are those among $t_1$, \dots, $t_n$ which are terms, plus the
immediate subterms of those among $t_1$, \ldots, $t_n$ which are not
terms.  Its sub-(semi)terms are the sub-(semi)terms of $t_1$, \ldots,
$t_n$.
\item If $A$ and $B$ are (semi)formulas, then $\lnot A$, $A \land B$, $A \lor
  B$ and $A \impl B$ are (semi)\-formulas. Its (immediate) sub-(semi)terms are
  those of $A$ and~$B$.
\item If $A(a)$ is a (semi)formula containing the free variable $a$ and $x$ is
  a bound variable not occurring in $A(a)$, then $\forall x\, A(x)$ and
  $\exists x\, A(x)$ are (semi)formulas. Its (immediate) sub-(semi)terms are
  those of~$A(x)$.
\item If $A(a)$ is a (semi)formula containing the free variable $a$ and $x$ is
  a bound variable not occurring in $A(a)$, then $\eps{x}{A(x)}$ is a
  (semi)term. Its sub-(semi)terms are $\eps{x}{A(x)}$ and the sub-(semi)terms
  of~$A(x)$. Its immediate sub-(semi)terms are those of~$A(x)$.
\end{enumerate}
Note that terms and formulas are just semiterm and semiformula which contain
no bound variables without a matching $\forall$, $\exists$, or $\meps$.  We
will call a semi-term of the form $\eps{x}{A(x)}$ an
\emph{$\meps$-expression}, or, if it is a term, an \emph{$\meps$-term}.
\end{definition}

The definition above does not allow a quantifier or epsilon binding a
variable~$x$ to occur in the scope of another quantifier or epsilon which
binds the same variable.  In order for some substitutions of terms to result
in well-formed formulas, it will often be necessary to rename bound variables.
Conversely, we cannot assume that all epsilon terms of the same form occurring
in a proof are syntactically identical.  Hence we will (a) adopt the
convention that bound variables must be renamed when substituting terms so as
to avoid clashes of bound variables, and (b) we tacitly identify epsilon terms
which differ only by a renaming of bound variables.  Thus, if $A(a)$ is a
formula with free variable $a$, then $A(\eps{x}{A(x)})$ is a formula obtained
from $A(a)$ by replacing every indicated occurrence of $a$ by an epsilon term
resulting from $\eps{x}{A(x)}$ by renaming each bound variable $y$ occurring
in it so as to ensure that $y$ does not appear in the scope of a quantifier or
epsilon in $A(a)$.  For instance, if $A(a) \equiv \exists y\, P(a, y)$, then
$\eps{x}{A(x)} \equiv \eps{x}{\exists y\, P(x, y)}$ contains the bound
variable $y$, and a literal substitution for $a$ in $A(a)$ would result in
$\exists y\, P(\eps{x}{\exists y\, P(x, y)}, y)$, which is not well formed.
The variable $y$ must be renamed, e.g., thus: $\exists y\, P(\eps{x}{\exists
  z\, P(x, z)}, y)$.

An instruction to replace every occurrence of an epsilon-term~$e$ in a formula
or proof by some other term is to be understood as an instruction to replace
every $\meps$-term equal to $e$ up to renaming of bound variables; e.g.,
``replace $\eps{x}{\exists y\, P(x, y)}$ in the above formula by $t$'' results
in $\exists y\, P(t, y)$.  It will sometimes be useful to be explicit about
substitutions of terms for variables.  We will write $A(a)\, \{a \gets t\}$ to
denote the result of replacing the indicated occurrences of $a$ in $A(a)$
by~$t$.

\Section{Axiomatization of the Epsilon Calculus} \label{sec:axioms}

\begin{definition}
  The set of axioms $\AxEC(L)$ of axioms of the \emph{elementary calculus} for
  a language $L$ consists of all propositional tautologies in the
  language~$L$. To obtain the set of axioms $\AxEC_{\meps}$ of the \emph{pure
    epsilon calculus} we add to $\AxEC(L(\ECeps))$ all substitution instances
  of
  \begin{equation} \label{eq:criticalformula}
    A(t) \impl A(\eps{x}{A(x)}) \tpkt
  \end{equation}
  An axiom of the form~\eqref{eq:criticalformula} is called a \emph{critical
    formula}.  We say that the critical formula \emph{belongs} to the
  $\meps$-term $\eps{x}{A(x)}$.
  
  The set $\AxPC$ of axioms of the \emph{predicate calculus}, and the set
  $\AxPC_{\meps}$ of axioms of the \emph{extended predicate calculus} consist
  of $\AxEC(L(\PC))$ and $\AxEC_\meps$, respectively, together with all
  instances of $A(t) \impl \exists x\, A(x)$ and $\forall x\, A(x) \impl A(t)$
  in the respective language.
\end{definition}

\begin{definition}
  A \emph{proof} in $\EC$ ($\ECeps$) is a sequence $A_1$, \dots, $A_n$ of
  formulas such that each $A_i$ is either in $\AxEC$ ($\AxEC_\meps$) or it
  follows from formulas preceding it by modus ponens, i.e., there are $j$, $k
  < i$ so that $A_k \equiv A_j \impl A_i$.
  
  A proof in $\PC$ ($\PCeps$) is a sequence $A_1$, \dots, $A_n$ of formulas
  such that each $A_i$ is either in $\AxPC$ ($\AxPC_\meps$) or follows from
  formulas preceding it by modus ponens, or follows from a preceding formula
  by generalization, i.e., there is a $j < i$ so that either $A_j \equiv B
  \impl C(a)$ and $A_i \equiv B \impl \forall x\, C(x)$ or $A_j \equiv B(a)
  \impl C$ and $A_i \equiv \exists x\, B(x) \impl C$.  In the latter case we
  also require that the free variable $a$ does not occur in $A_i$ or in any
  $A_k$ with $k > i$. Such a variable $a$ is called an \emph{eigenvariable}.
  The restriction guarantees that each variable is only used as an
  eigenvariable in a generalization inference at most once.
  
  A formula~$A$ is called \emph{provable} in $\EC$ ($\ECeps$, $\PC$, $\PCeps$)
  if there is a proof in $EC$ ($ECeps$, $\PC$, $\PCeps$, respectively) which
  has~$A$ as its last formula.  To indicate that $A$ is provable in, say,
  $\ECeps$ by a proof $\pi$, we write $\ECeps \vdash_{\pi} A$.
\end{definition}

\begin{definition}
  The \emph{size}~$\size{\pi}$ of a proof $\pi$ is the number of steps
  in~$\pi$.  If $\pi$ is a proof in $\ECeps$ or $\PCeps$, we define the
  \emph{critical count}~$\cc{\pi}$ of $\pi$ as the number of distinct critical
  formulas and quantifier axioms in~$\pi$ plus~$1$.
\end{definition}

\Section{The Embedding Lemma} \label{s:emblemma}

The epsilon operator allows the treatment of quantifiers in a quantifier-free
system: using $\meps$-terms, it is possible to define $\exists x$ and
$\forall x$ as follows:
\begin{eqnarray*}
\exists x\, A(x) & \Leftrightarrow & A(\eps{x}{A(x)}) \\
\forall x\, A(x) & \Leftrightarrow & A(\eps{x}{\lnot A(x)})
\end{eqnarray*}
We define a mapping $^\meps$ of semiformulas and semiterms in
$L(\PC_\meps)$ to semiformulas and semiterms in $L(\ECeps)$ as
follows:
\[
\begin{array}{@{}rclrclrcl@{}}
\multicolumn{9}{c}{f(t_1, \ldots, t_n)^\meps =
    f(t_1^\meps, \ldots, t_n^\meps)
 \qquad P(t_1,\ldots, t_n)^\meps  = P(t_1^\meps, \ldots, t_n^\meps)} \\
x^\meps & = & x 
& (A \impl B)^\meps & = & A^\meps \impl B^\meps 
& [\eps{x}{A(x)}]^\meps & = & \eps{x}{A(x)^\meps}  \\
 a^\meps & = & a & 
(A \lor B)^\meps & = & A^\meps  \lor B^\meps & 
(\exists x\, A(x))^\meps & = & A^\meps(\eps{x}{A(x)^\meps}) \\
(\lnot A)^\meps  & = & \lnot A^\meps &
(A \land B)^\meps & = & A^\meps \land B^\meps & 
(\forall x\, A(x))^\meps & = & A^\meps(\eps{x}{\lnot A(x)^\meps})
\end{array}
 \]
\begin{example} Consider
\begin{eqnarray*}
\exists x(P(x) & \lor & \forall y\, Q(y))^\meps = \\ 
& = & [P(x) \lor \forall y\, Q(y)]^\meps \quad 
\{x \gets \eps{x}{[P(x) \lor \forall y\, Q(y)]^\meps}\} \\ 
&& \quad [P(x) \lor \forall y\, Q(y)]^\meps = 
   P(x) \lor Q(\underbrace{\eps{y}{\lnot Q(y)}}_{e_1}) \\ 
& = & P(x) \lor Q(\underbrace{\eps{y}{\lnot Q(y)}}_{e_1}) \quad
  \{x \gets \underbrace{\eps{x}{[P(x) \lor 
        Q(\underbrace{\eps{y}{\lnot Q(y)}}_{e_1})]}}_{e_2}\} \\ 
& = & P(\underbrace{\eps{x}{[P(x) \lor 
      Q(\underbrace{\eps{y}{\lnot Q(y)}}_{e_1})]}}_{e_2}) 
         \lor Q(\underbrace{\eps{y}{\lnot Q(y)}}_{e_1})
\end{eqnarray*}
\end{example}

\begin{example} Consider
\begin{eqnarray*}
[\exists x\, & \exists y  & A(x, y)]^\meps = \\
 & = & [\exists y\, A(x, y)]^\meps \quad
\{x \gets \eps{x}{[\exists y\, A(x, y)]^\meps}\} \\
& & \quad [\exists y\, A(x, y)]^\meps = 
   A(x, \underbrace{\eps{y}{A(x, y)}}_{e'(x)}) \\
& = & A(x, \underbrace{\eps{y}{A(x, y)}}_{e'(x)}) 
  \{x \gets \underbrace{\eps{x}{A(x, \eps{z}{A(x, z)})}}_{e_3}\} \\
& = & A(\underbrace{\eps{x}{A(x, \eps{z}{A(x, z)})}}_{e_3}, 
          \underbrace{\eps{y}{A(\underbrace{\eps{x}{A(x, 
                  \eps{z}{A(x, z)})}}_{e_3}, y)}}_{e_4\,=\,e'(e_3)})
\end{eqnarray*}  
\end{example}

\begin{lemma}[Embedding Lemma]
  If $\pi$ is a $\PCeps$-proof of $A$ then there is an $\ECeps$-proof
  $\pi^\meps$ of $A^\meps$ with $\size{\pi^\meps} \le 3\cdot \size{\pi}$ and
  $\cc{\pi^\meps} \le \cc{\pi}$.
\end{lemma}

\begin{proof}
  We show that for all proofs $\pi$ consisting of formulas $A_1$, \ldots,
  $A_n$, there is a proof $\pi^\meps$ containing $A_1^\meps$, \ldots,
  $A_n^\meps$ (plus perhaps some extra formulas) of the required size and
  critical count.  We proceed by induction on~$n$.  The case $n = 0$ is
  trivial.  Suppose the claim holds for the proof consisting of $A_1$, \ldots,
  $A_{n}$, i.e., there is a proof $\pi^*$ containing $A_1^\meps$, \ldots,
  $A_{n}^\meps$, and consider the proof $\pi = A_1$, \dots, $A_n$, $A$.  If
  $A$ is a propositional tautology, then $A^\meps$ is also a propositional
  tautology.  (Note that $(\cdot)^\meps$ leaves the propositional structure of
  $A$ intact.) If $A$ is a critical formula, then $A^\meps$ is also a critical
  formula.  In both cases, we can take $\pi^\meps$ to be $\pi^*$ extended by
  $A_n^\meps$.
  
  If $A$ is an instance of a quantifier axiom, its translation $A^\meps$
  either is of the form
\begin{equation*}
[A(t) \impl \exists x\, A(x)]^\meps \equiv 
A^\meps(t^\meps) \impl A^\meps(\eps{x}{A(x)^\meps}) \tkom
\end{equation*}
which is a critical formula, or is of the form
\begin{equation*}
[\forall x\, A(x) \impl A(t)]^\meps \equiv
   A^\meps(\eps{x}{\lnot A(x)}) \impl A^\meps(t^\meps) \tkom
\end{equation*}
which is the contrapositive of, and hence provable from, a critical formula.
In the latter case, the size of the resulting proof increases by two
additional steps.

Now suppose $A$ follows by modus ponens from $A_i$ and $A_j \equiv A_i \impl
A$.  Since $\pi^*$ contains $A_i^\meps$ and $A_j^\meps \equiv A_i^\meps \impl
A^\meps$, adding $A^\meps$ to $\pi^*$ is also a proof.

If $A$ follows by generalization, i.e., $A \equiv B \impl \forall x\, C(x)$
and $A_i \equiv B \impl C(a)$ (where $a$ satisfies the conditions on
eigenvariables), then by induction hypothesis the proof $\pi^*$ contains
$A_i^\meps \equiv B^\meps \impl C(a)^\meps$. Replacing $a$ everywhere in
$\pi^*$ by $\eps{x}{\neg A(x)}$ results in a proof containing
\begin{equation*}
[B \impl \forall x\, C(x)]^\meps \equiv B^\meps \impl A^\meps(\eps{x}{\lnot
A(x)^\meps}) \tpkt
\end{equation*}
in place of $A_i^meps$. Similarly, if the last inference derives $A \equiv
\exists x\, B(x) \impl C$ from $B(a) \impl C$, by induction hypothesis $\pi^*$
contains $B(a)^\meps \impl C^\meps$, and we obtain a proof of $A^\meps$ by
replacing $a$ everywhere in~$\pi^*$ by $\eps{x}{B(x)}$.
\end{proof}

\Section{The First Epsilon Theorem}

We begin our discussion of the $\meps$-theorems by a detailed proof of the
first $\meps$-theorem.  This theorem states that if a formula~$E$ without
quantifiers or epsilon is provable in the (extended) $\meps$-calculus, then it
is already provable in the elementary calculus.  In other words, the
(extended) epsilon calculus is conservative over the elementary calculus for
elementary formulas.  A relatively simple corollary of the first epsilon
theorem is the extended first $\meps$-theorem, which is a version of
Herbrand's theorem for prenex formulas.  This section is dedicated to proofs
of these results.  The argument we use is essentially Bernays's~\cite{HB70},
which gives a procedure by which critical formulas in proofs of~$E$ are
eliminated step-wise.  We analyze the procedure and thereby obtain upper
bounds on the complexity of the Herbrand disjunction obtained in the extended
first $\meps$-theorem.  These bounds are given in terms of the
\emph{hyperexponential function}~$2_y^x$, defined by $2_0^x = x$ and
$2_{i+1}^x = 2^{2_i^x}$.

The proof will proceed by induction on the rank and degree and number of
$\epsilon$-terms in critical formulas in the proof of~$E$.  Rank and degree
are two measures of complexity of $\meps$-expressions: Degree applies to
$\meps$-terms only, and measures the depth of nesting of $\meps$-terms.  Rank,
on the other hand, measures the complexity of cross-binding of
$\meps$-expressions.

\begin{definition}
  The \emph{degree} of an $\meps$-term is inductively defined as follows:
  \begin{enumerate}
  \item If $A(x)$ contains no $\meps$-subterms, then $\deg(\eps{x}{A(x)}) =
    1$.
  \item If $e_1$, \dots, $e_n$ are all immediate $\meps$-subterms of $A(x)$,
    then
    \begin{equation*}
      \deg(\eps{x}{A(x)}) = \max\{\deg(e_1), \ldots,
      \deg(e_n)\} + 1 \tpkt  
    \end{equation*}
  \end{enumerate}
\end{definition}

\begin{definition}
  An $\meps$-expression $e$ is \emph{subordinate} to $\eps{x}{A(x)}$ if $e$ is
  a proper sub-semiterm of $A(x)$ and contains~$x$.
\end{definition}

\begin{definition}
  The \emph{rank} of an $\meps$-expression $e$ is defined as follows:
  \begin{enumerate}
  \item If $e$ contains no subordinate $\meps$-expressions, then $\rk(e)
    = 1$.
  \item If $e_1$, \ldots, $e_n$ are all the $\meps$-expressions subordinate to
    $e$, then
    \begin{equation*}
      \rk(e) = \max\{\rk(e_1), \ldots, \rk(e_n)\} + 1 \tpkt 
    \end{equation*}
  \end{enumerate}
\end{definition}

\begin{example}
First, consider the formula
\begin{equation*}
P(\underbrace{\eps{x}{[P(x) \lor 
    Q(\underbrace{\eps{y}{\lnot Q(y)}}_{e_1})]}}_{e_2}) \lor 
    Q(\underbrace{\eps{y}{\lnot Q(y)}}_{e_1}) \tpkt
\end{equation*}
Here, $e_1$ is the only immediate $\meps$-subterm of $e_2$ and has no
$\meps$-subterms itself, so $\deg(e_1) = 1$ and $\deg(e_2) = 2$.  Neither
$e_1$ nor $e_2$ contains subordinate $\meps$-expressions, hence $\rk(e_1) =
\rk(e_2) = 1$.  In
\begin{equation*}
[\exists x\exists y\, A(x, y)]^\meps = A(\underbrace{\eps{x}{A(x, \eps{z}{A(x,
      z)})}}_{e_3}, \underbrace{\eps{y}{A(\underbrace{\eps{x}{A(x,
        \eps{z}{A(x, z)})}}_{e_3}, y)}}_{e_4}) \tkom
\end{equation*}
$e_3$ contains no $\meps$-subterms, but $e_4$ contains $e_3$ as a subterm, so
$\deg(e_3) = 1$ and $\deg(e_4) = 2$.  On the other hand, $e_3$ contains the
subordinate $\meps$-expression $\eps{z}{A(x, z)}$, hence $\rk(e_3) = 2$.
Since $y$ does not occur in the scope of another $\meps$, $e_4$ contains no
subordinate $\meps$-expressions, and $\rk(e_4) = 1$.
\end{example}

\begin{definition}
  Suppose $\pi$ is a proof in $\ECeps$.  If $e$ is an $\meps$-term belonging
  to a critical formula $A(t) \impl A(e)$ of $\pi$, then we call $e$ a
  \emph{critical epsilon term} of $\pi$, and $\rk(e)$ ($\deg(e)$) the rank
  (the degree) of that critical formula.
  
  The \emph{rank}~$\rk(\pi)$ of $\pi$ is the maximum rank of its critical
  formulas.  The \emph{degree} $\deg(\pi, r)$ of~$\pi$ with respect to
  rank~$r$ is the maximum degree of its critical $\meps$-terms of rank~$r$.
  The \emph{order} $o(\pi,r)$ of $\pi$ with respect to rank~$r$ is the number
  of different critical $\meps$-terms of rank~$r$.
\end{definition}

\begin{lemma} \label{l:a}
  Let $\pi$ be a $\ECeps$-proof, let $r = \rk(\pi)$ be the maximal rank
  of critical formulas in $\pi$, and let $e$ be a critical $\meps$-term of
  $\pi$ of maximal degree among the critical $\meps$-terms of rank~$r$.
  
  Suppose that $A(t) \impl A(e)$ is a critical formula belonging to~$e$ and
  that $B^* \equiv B(s) \impl B(\eps{y}{B(y)})$ is a critical formula in $\pi$
  belonging to a different $\meps$-term~$\eps{y}{B(y)}$, and suppose $C$ is
  the result of replacing $e$ by $t$ in~$B^*$. Then (a) if $\rk(B^*) = r$,
  then $C$ and $B^*$ have the the same critical $\meps$-term $\eps{y}{B(y)}$
  belonging to them, and (b) $\rk(C) = \rk(B^*)$.
\end{lemma}

\begin{proof}
  We first consider and exclude a preliminary case. If the indicated
  occurrences of $s$ (on the left-hand side) or the indicated occurrences of
  $\eps{y}{B(y)}$ (on the right-hand side) lie \emph{inside}~$e$, replacing
  $e$ by $t$ would result in a formula which is not of the form of a critical
  formula. This, however, can never be the case.  For suppose it were, i.e.,
  suppose $B(a)$ is of the form $B'(e'(a))$ and either $e \equiv e'(s)$ or $e
  \equiv e'(\eps{y}{B(y)})$.  If $e \equiv e'(s)$, the left-hand-side $B(s)$
  of the critical formula is of the form $B'(e'(s))$, and consequently the
  right-hand-side is $B'(e'(\eps{y}{B'(e'(y))}))$.  Conversely, if $e \equiv
  e'(\eps{y}{B(y)})$, then the right-hand-side $B(\eps{y}{B(y)})$ would be of
  the form $B'(e'(\eps{y}{B'(e'(y))}))$.  In either case we have an
  $\meps$-term~$e'(a)$ which is of the same rank as $e$, since $e = e'(t')$
  for some term~$t'$. (Note that by our conventions on renaming of variables,
  no subexpression of $t'$ can be subordinate to~$e$.)  On the other hand,
  $e'(y)$ is subordinate to the critical $\meps$-term $\eps{y}{B(y)} =
  \eps{y}{B'(e'(y))}$, and so $\rk(\eps{y}{B(y)}) > \rk(e'(y)) = \rk(e)$. But
  $e$ was assumed to be a critical $\meps$-term of maximal rank.
  
  There are then only two ways in which $e$ can occur in a critical formula:
  either (i) $e$ occurs only in the indicated occurrences of~$s$ but not in
  $B(y)$ at all, or (ii) $e$ occurs in $B(y)$ (and perhaps also in $s$). 
  
  Case (i): $e$ occurs only in $s$, i.e., $s \equiv s'(e)$.  Replacing $e$ by
  $t$ results in the critical formula $C \equiv B(s'(t)) \impl
  B(\eps{y}{B(y)})$.  The new critical formula $C$ belongs to the same
  $\meps$-term as the original formula, hence we obtain (a) $\rk(C) = \rk(B^*)
  = \rk(\eps{y}{B(y)})$, and (b) holds trivially.

  Case~(ii): $e$ occurs in $B(y)$ (and perhaps also in $s$). In this case,
  $B(y) \equiv B'(y, e)$ and the critical formula has the form
  \begin{equation*}
    B'(s, e) \impl B'(\eps{y}{B'(y, e)}, e) \tpkt
  \end{equation*}  
  Then the $\meps$-term belonging to this critical formula, $e' \equiv
  \eps{y}{B'(y, e)}$, contains $e$ as a proper subterm and hence is of higher
  degree than~$e$.  Since $e$ is a critical $\meps$-term of maximal degree
  among the critical $\meps$-terms of maximal rank in $\pi$, this implies that
  $\rk(e') < \rk(e)$.  Replacing $e$ by $t$ yields the critical formula
  \begin{equation*}
    C \equiv B'(s', t) \impl B'(\eps{y}{B'(y, t)}, t) \tkom  
  \end{equation*}
  belonging to the $\meps$-term $\eps{y}{B'(y, t)}$.  This term has the same
  rank as $e'$ and hence a lower rank than $e$ itself (although it might have
  a degree higher than~$\deg(e)$).  Hence, $\rk(C) < r$.
\end{proof}

\begin{lemma} \label{l:1}
  Suppose $\ECeps \vdash_\pi E$. Let $r = \rk(\pi)$ be the maximal rank of
  critical formulas in $\pi$, and let $e$ be a critical $\meps$-term of $\pi$
  of maximal degree among the critical $\meps$-terms of rank~$r$.  Then there
  is a proof $\pi_e$ so that $\ECeps \vdash_{\pi_e} E$ with $\rk(\pi_e) \le
  r$, $\deg(\pi_e, r) \le \deg(\pi, r)$ and $o(\pi_e, r) = o(\pi, r) - 1$.
\end{lemma}

\begin{proof}
  The $\meps$-expression $e$ is of the form $\eps{x}{A(x)}$, and suppose that
  $A(t_1) \impl A(e)$, \dots, $A(t_n) \impl A(e)$ are all the critical
  formulas in $\pi$ belonging to~$e$.  For each $i = 1$, \dots, $n$, we obtain
  a proof $\pi_i$ of $A(t_i) \impl E$ as follows:
  \begin{enumerate}
  \item Replace $e$ everywhere in $\pi$ by $t_i$.  Every critical formula
    $A(t_j) \impl A(e)$ belonging to $e$ thus turns into a formula of the form
    $A(t_j') \impl A(t_i)$.  To see this, note that $e$ cannot occur in $A(x)$,
    for otherwise $e \equiv \eps{x}{A(x)}$ would be a proper subterm of itself,
    which is impossible.
  \item Add $A(t_i)$ to the axioms. Now every one of the new formulas $A(t_j')
    \impl A(t_i)$ is derivable using modus ponens from the tautology 
    \begin{equation*}
      A(t_i) \impl (A(t_j') \impl A(t_i)) \tkom
    \end{equation*}
  \item Apply the deduction theorem to obtain~$\pi_i$.
  \end{enumerate}
  We verify that $\pi_i$ is indeed an $\ECeps$-proof with the required
  properties.  In the construction of $\pi_i$, we replaced $e$ by $t_i$
  throughout the proof.  Such a substitution obviously preserves tautologies.
  By Lemma~\ref{l:a}, it also turns critical formulas into critical formulas
  (belonging, perhaps, to different critical $\meps$-terms).  Lemma~\ref{l:a}
  also guarantees that replacing $e$ by $t_i$ everywhere does not change the
  rank of critical formulas, and that it does not change the critical
  $\meps$-terms of maximal rank~$r$ at all (in particular, it does not
  increase the degree of critical $\meps$-terms of rank~$r$).
  
  We started with critical formulas $A(t_i) \impl A(e)$, and obtained a
  proof~$\pi_i$ which does not contain any critical formulas belonging to~$e$.
  Hence $e$ is no longer a \emph{critical} $\meps$-term in $\pi_i$. The ranks
  of all other critical formulas (and the corresponding critical
  $\meps$-terms) remain unchanged. Thus $o(\pi_i, r) = o(\pi, r)-1$.
  
  Secondly, we construct a proof $\pi'$ of $\bigwedge_{i=1}^n \lnot A(t_i)
  \impl E$ as follows:
  \begin{enumerate}
  \item Add $\bigwedge \lnot A(t_i)$ to the axioms. Now every critical formula
    $A(t_i) \impl A(e)$ belonging to $e$ is provable using the propositional
    tautology $\lnot A(t_i) \impl (A(t_i) \impl A(e))$.
  \item Apply the deduction theorem for the propositional calculus to obtain a
    proof~$\pi'$, which contains exactly the same critical formulas as~$\pi$
    except those belonging to~$e$.
  \end{enumerate}
  Now combine the proofs $\pi_i$ of $A(t_i) \impl E$ and $\pi'$ of
  $\bigwedge_i \lnot A(t_i) \impl E$ to get the proof~$\pi_e$ of $E$ (using
  case distinction).
  
  None of the proofs $\pi_i$, $\pi'$ contain critical formulas belonging
  to~$e$, and no critical $\meps$-terms of rank~$r$ other than those in~$\pi$.
  Thus $\rk(\pi_e) \le r$, $\deg(\pi_e, r) \le \deg(\pi, r)$, and $o(\pi_e,r)
  = o(\pi,r)-1$ hold.
\end{proof}

\begin{theorem}[First Epsilon Theorem] \label{t:1} 
  If $E$ is a formula without bound variables (no quantifiers, no epsilons)
  and $\PCeps \vdash E$ then $\EC \vdash E$.
\end{theorem}

\begin{proof}
  First, use the embedding lemma to obtain $\pi$ so that $\ECeps \vdash_\pi
  A$.  The theorem then follows from the preceding lemma by induction on $r =
  \rk(\pi)$ and $d = o(\pi, r)$.  If $r = 0$, there are no critical formulas,
  so there is nothing to prove.  If $r > 0$, then $d$-fold application of the
  lemma results in a proof~$\pi'$ of rank~$< r$.
\end{proof}

\begin{theorem}[Extended First Epsilon Theorem]\label{t:1e}
  Suppose $E(e_1, \ldots, e_m)$ is a quantifier-free formula containing only
  the $\meps$-terms $s_1$, \ldots, $s_m$, and
  \begin{equation*}
    \ECeps \vdash_\pi E(s_1, \ldots, s_m) \tkom
  \end{equation*}
  then there are $\meps$-free terms 
  $t^i_j$ ($1 \le i  \le n$, $1 \le j \le m$) such that
  \begin{equation*}
    \EC \vdash \bigvee_{i=1}^n E(t^i_1, \ldots, t^i_m)
  \end{equation*}
  where $n \leq 2_{2\cdot\cc{\pi}}^{3\cdot\cc{\pi}}$.
\end{theorem}

\begin{corollary}[Herbrand's Theorem]\label{c:herbrand1} 
  If $\exists x_1\ldots \exists x_m E(x_1, \ldots, x_m)$ is a purely
  existential formula containing only the bound variables $x_1$, \dots, $x_m$,
  and
  \begin{equation*}
    \PCeps \vdash_\pi \exists x_1\ldots \exists x_m E(x_1, \ldots, x_m) \tkom
  \end{equation*}
  then there are $\meps$-free terms 
  $t^i_j$ ($1 \le i  \le n$, $1 \le j \le m$) such that
  \begin{equation*}
    \EC \vdash \bigvee_{i=1}^n E(t^i_1, \ldots, t^i_m)
  \end{equation*}
  where $n \leq 2_{2\cdot\cc{\pi}}^{3\cdot\cc{\pi}}$.
\end{corollary}

\begin{proof} Immediate from Theorem~\ref{t:1e} using the Embedding
  Lemma. \end{proof}

The rest of this section is devoted to the proof of 
of Theorem~\ref{t:1e}. First, some additional notation.  

\begin{definition}
  Suppose $\pi$ is a proof in $\ECeps$. The \emph{width}~$\width_{\pi}(e)$ of
  $\pi$ with respect to~$e$ is the number of different critical formulas in
  $\pi$ belonging to the $\meps$-term~$e$. The \emph{width}~$\width(\pi, r)$ of
  $\pi$ with respect to rank~$r$ is given by
\begin{equation*}
    \width(\pi,r) =
     \max \{\width_{\pi}(e) \mid 
           \text{$e$ of rank $r$ occurs in $\pi$} \} + 1\tpkt
\end{equation*}
\end{definition}

\begin{definition}
Let $E(a_1, \ldots, a_m)$ be a formula in $L(\EC)$ without bound variables,
and let $s_1$, \ldots, $s_m$ be terms in $L(\ECeps$). An
\emph{$\lor$-expansion} (of $E \equiv E(s_1,\ldots,s_m)$) is a finite
disjunction
\begin{equation*}
    E'\equiv  E_1 \lor \cdots \lor E_l \tkom
\end{equation*}
where each $E_i \equiv E(s^i_1,\ldots,s^i_m)$ for terms $s^i_j$ ($1 \le i
\le l$, $1 \le j \le m$).  We call $l$ the \emph{length} $\expan(E',E)$ of the
expansion.
\end{definition}

\begin{proposition}\label{p:expanlen}
  Suppose $A'$ is an $\lor$-expansion of $A$ and $A''$ is an $\lor$-expansion
  of $A'$.  Then $A''$ is also an $\lor$-expansion of $A$, and $\expan(A'',A)
  \leq \expan(A'',A') \cdot \expan(A',A)$.
\end{proposition}

\begin{proof} Obvious. \end{proof}

\begin{lemma} \label{l:2}
  Suppose $E(a_1, \ldots, a_m)$ is a formula in $L(\EC)$ without bound
  variables and $\pi$ is an $\ECeps$-proof of $E(s_1, \ldots, s_m)$ where
  $s_1$, \ldots, $s_m$ are terms in $L(\ECeps)$. Let $r = \rk(\pi)$ and let
  $e$ be a critical $\meps$-term of $\pi$ of maximal
  degree among the critical $\meps$-terms of rank~$r$, and let $n =
  \width_{\pi}(e)$ be the number of critical formulas belonging to~$e$.  Then
  there are terms $s^i_j$ $(0 \le i \le n+1$, $1 \le j \le m)$ and a proof
  $\pi_e$ with end formula
  \begin{equation*}
    E(s^1_1, \ldots, s^1_m) \lor \cdots \lor 
    E(s^{n+1}_1, \ldots, s^{n+1}_m) \tkom
  \end{equation*}
  so that $\rk(\pi_e) \le r$, $\deg(\pi_e, r) \le \deg(\pi, r)$, and $o(\pi_e,
  r) = o(\pi, r) - 1$.  Furthermore $\cc{\pi_e} \leq \cc{\pi} \cdot (n+1)$ and
  $\width(\pi_e, r) \leq \width(\pi,r) \cdot (n+1) \leq \width(\pi,r)^2$.
\end{lemma}

\begin{proof}
  We adapt the construction of $\pi_i$ in Lemma~\ref{l:1}.  The only
  difference to the previous construction is that when replacing $e$ by $t_i$
  throughout~$\pi$, the end-formula $E(s_1, \ldots, s_m)$ may change.
  However, $e$ can only occur in $s_1$,~\ldots,~$s_m$ since $E(a_1, \ldots,
  a_m)$ contains no bound variables and hence no $\meps$-terms.  For each
  critical formula $A(t_i) \impl A(e)$ we obtain a proof $\pi_i$ of $A(t_i)
  \impl E(s^i_1, \ldots, s^i_m)$. The construction of $\pi'$ as before yields
  a proof of
  \begin{equation*}
    \bigwedge_{i=1}^n \lnot A(t_i) \impl E(s^{n+1}_1, \ldots, s^{n+1}_m) \tkom
  \end{equation*}
  if we take $s_j^{n+1} = s_j$.  Then, since obviously for each~$i$
  \begin{equation*}
    E(s^i_1, \ldots, s^i_m) \impl \bigvee_{i=1}^{n+1} E(s^i_1, 
  \ldots, s^i_m) \tkom
  \end{equation*}
  is provable, we obtain a proof~$\pi_e$ of $\bigvee_{i=1}^{n+1} E(s^i_1,
  \ldots, s^i_m)$ with the desired properties.  Observe that the length of the
  $\lor$-expansion $\bigvee E(s^i_1, \ldots, s^i_m)$ is $n+1 =
  \width_{\pi}(e)+1 \leq \width(\pi,r)$.
  
  It remains to verify the bounds on $\cc{\pi_e}$ and $\width(\pi_e,r)$.  By
  Lemma~\ref{l:a}, replacing $e$ by $t_i$ to obtain $\pi_i$ does not introduce
  new critical $\meps$-terms of rank~$r$.  (Critical formulas belonging to
  $\meps$-terms of $\rk(e)$ may be altered, but the corresponding critical
  $\meps$-terms remain the same.)  New critical $\meps$-terms can only appear
  at a rank $< \rk(e)$, and if they do, their rank is equal to the rank of a
  some critical $\meps$-term already in~$\pi$.  The total number of critical
  formulas in $\pi_i$ is at most that of $\pi$ less the number~$n$ of critical
  formulas belonging to~$e$, i.e., $\cc{\pi_i} \leq \cc{\pi}-n$.  Moreover,
  obviously $\cc{\pi'} \leq \cc{\pi}-n$ holds.
  
  When we combine the $n+1$ proofs $\pi_i$ and $\pi'$ to obtain $\pi_e$, the
  worst case is that every critical formula in~$\pi_i$ has been changed. Thus
  $\cc{\pi_e} \le (\cc{\pi}-n)(n+1) \le \cc{\pi}(n+1)$.  Now looking more
  closely at the critical formulas of rank~$r$ in $\pi_i$, we see that
  whenever case (i) in the proof of Lemma~\ref{l:a} applies, a critical
  formula belonging to some $\meps$-term~$e'$ of rank~$r$ in $\pi$ turns into
  a potentially new critical formula in~$\pi_i$ also belonging to~$e'$.
  However, these are the only new critical formulas belonging to~$e'$.  Hence,
  there may be up to $\width_\pi(e')\cdot(n+1)$ different critical formulas
  belonging to~$e'$ in $\pi_e$.  Thus $\width(\pi_e,r) \le
  \width(\pi,r)\cdot(n+1)$, which is $\le \width(\pi,r)^2$ since $n+1 =
  \width_\pi(e) +1 \le \width(\pi,r)$.
\end{proof}

We now iterate the elimination of $\meps$-terms of highest rank and estimate
the critical count of the proof resulting from the elimination of all
$\meps$-terms of rank $\rk(e)$.

\begin{lemma} \label{l:b2}
  Suppose $E(a_1, \ldots, a_m)$ is a formula in $L(\EC)$ without bound
  variables and $\pi$ is an $\ECeps$-proof of $E(s_1, \ldots, s_m)$, where
  $s_1$, \ldots, $s_m$ are terms in $L(\ECeps$).  Then there is a proof
  $\sigma$ of an $\lor$-expansion $E'$ of $E(s_1, \ldots, s_m)$, so that
  $\rk(\sigma) < \rk(\pi)$.  Furthermore,
  \begin{equation*}
    \cc{\sigma} \leq 2^{2^{2\cdot \cc{\pi}}} \qquad \text{and} \qquad
    \expan(E',E(s_1, \ldots, s_m)) \leq 2^{2^{2\cdot \cc{\pi}}} \tpkt
  \end{equation*}
\end{lemma}

\begin{proof}
  Let $d = o(\pi,r)$ and let $e_1$, \ldots, $e_d$ be all critical
  $\meps$-terms of rank~$r$ in~$\pi$.  We assume the sequence $e_1$, \dots,
  $e_d$ is ordered so that the degree never increases.  Let $k = \cc{\pi}$
  and $r = \rk(\pi)$.  As observed in the preceding proof, an application of
  Lemma~\ref{l:2} cannot increase the number of critical $\meps$-terms of
  maximal rank.  Thus let $\sigma^0 = \pi$ and $\sigma^{j} =
  \sigma^{j-1}_{e_j}$ for $j>0$.  We thus obtain $\sigma = \sigma^d$ by
  $d$-fold iteration of Lemma~\ref{l:2}.  Let $E^0 \equiv E(s_1, \ldots,
  s_m)$. By construction, the critical $\meps$-terms of rank~$r$ in $\sigma^j$
  are $e_j$, \dots, $e_d$ and the end-formula~$E^j$ of $\sigma^j$ is an
  $\lor$-expansion of $E \equiv E(s_1, \ldots, s_m)$; we set $E' \equiv E^d$.
  
  By induction on $j$ we prove:
  \begin{eqnarray*}
    \width(\sigma^j, r) & \leq & k^{2^j} \tkom \\
    \cc{\sigma^j} & \leq & k \cdot k^{\sum_{l=0}^{j-1} 2^l} \tkom \\
    \expan(E^j,E) & \leq & k^{\sum_{l=0}^{j-1} 2^l} \tpkt
  \end{eqnarray*}
  For $j=0$, we have $\width(\sigma^0,r) = \width(\pi, r) \le k$,
  $\cc{\sigma^0} = \cc{\pi} = k$ and $\expan(E^0, E(s_1, \ldots,s_m)) = 1$, by
  definition.

                       
  Now assume that $\sigma^j$ obeys the stated bounds, we prove that
  $\sigma^{j+1}$ does as well.  Apply Lemma~\ref{l:2} to
  eliminate the $\meps$-term $e_{j+1}$.  This yields a proof $\sigma^{j+1}$
  of $E^{j+1}$.  By Lemma~\ref{l:2} and the induction hypothesis we have:
  \begin{eqnarray*}
    \width(\sigma^{j+1}, r) & \leq & \width(\sigma^j, r)^2 \leq
    (k^{2^j})^2 = k^{2^{j+1}} \tkom \\
    \cc{\sigma^{j+1}} & \leq & (k \cdot k^{\sum_{l=0}^{j-1} 2^l}) \cdot
    k^{2^j} = k \cdot k^{\sum_{l=0}^{j} 2^l}, \text{and} \\
    \expan(E^{j+1}, E) & \leq & (k^{\sum_{l=0}^{j-1} 2^l}) \cdot 
    k^{2^j} = k^{\sum_{l=0}^{j} 2^l} \tpkt
  \end{eqnarray*}
  since $\width_{\sigma^{j}}(e_{j+1}) + 1 \le \width(\sigma^j, r)
  \le k^{2^j}$.
  Thus the claim follows and we obtain
  \begin{eqnarray*}
    \cc{\sigma} & \leq & k \cdot k^{\sum_{l=0}^{d-1} 2^l} =
    k \cdot k^{2^d - 1} = k^{2^d} \quad \text{and}\\
    \expan(E',E) & \leq & k^{\sum_{l=0}^{d-1} 2^l} = k^{2^d-1} \tpkt
  \end{eqnarray*}
  The order $d = o(\pi,r)$ is the number of critical $\meps$-terms of rank
  $r$, and hence $\le k$. Using the inequality $k \le 2^k$ we obtain the
  (rough) upper bounds $\cc{\sigma} \leq k^{2^d} \leq 2^{2^{2k}}$ and
  $\expan(E',E) \leq k^{2^{d-1}} \leq 2^{2^{2k}}$.
\end{proof}

\begin{proof}[Proof of Extended First Epsilon Theorem~\ref{t:1e}]
  Consider a proof $\pi$ of $E(s_1, \ldots, s_m)$, where $s_1$, \ldots, $s_m$
  are terms containing $\meps$'s and $E(a_1, \ldots, a_m)$ contains no bound
  variables. Let $k = \cc{\pi}$ and $p$ be the number of different ranks of
  critical $\meps$-terms in $\pi$, i.e., $p = \card{\{r : \width(\pi, r) >
    1\}}$.  The number $p$ is the number of times we have to apply
  Lemma~\ref{l:b2} to eliminate all critical formulas from $\pi$.  (Note that
  by the proof of Lemma~\ref{l:2} each elimination step can only decrease the
  number of different ranks.)  Although obviously $p \le r = \rk(\pi)$ we also
  have $p \le k = \cc{\pi}$, so the number of times Lemma~\ref{l:b2} must be
  applied is actually independent of~$\rk(\pi)$.
  
  Now let $\pi^0 = \pi$ and $\pi^{j+1}$ be the proof $\sigma$ constructed in
  Lemma~\ref{l:b2} starting with $\pi^j$.  Note that the end-formula of each
  $\pi^j$ is an $\lor$-expansion of~$E$. If we write $E^j$ for the
  end-formula of~$\pi^j$, then $E^p$ is the required Herbrand disjunction
  \begin{equation}\label{eq:p:h1}
    \bigvee_{i=1}^n E(t_1^i, \ldots, t_m^i) \tpkt
  \end{equation}
  To establish a bound on the length $n$ of~\eqref{eq:p:h1}, we apply
  Lemma~\ref{l:b2} $(p-1)$ times.  This yields the following bounds:
  \begin{equation*}
    \cc{\pi^{p-1}} \leq 2_{2(p-1)}^{2k+(p-1)} \quad \text{and} \quad
    \expan(E_{p-1},E) \leq 2_{2(p-1)}^{2k+(p-1)} \tpkt
  \end{equation*}
  Another application of Lemma~\ref{l:b2} yields that 
  \begin{eqnarray*}
    n & \leq & \expan(E_{p},E_{p-1}) \cdot \expan(E_{p-1},E) \\
     & \leq & 2^{2^{2\cdot \cc{\pi^{p-1}}}} \cdot \expan(E_{p-1},E) \\
     & \leq & 2^{2^{2\bigl(2_{2(r-1)}^{2k+(r-1)}\bigr)}} \cdot 2_{2(p-1)}^{2k+(p-1)} \\
     & \leq & 2^{2^{2_{2(p-1)}^{2k+p}}} = 
              2_{2p}^{2k + p} \leq 2_{2k}^{3k} \tpkt
  \end{eqnarray*}
  As a last step, we remove the remaining (non-critical) $\meps$-terms from
  the proof by replacing outermost $\meps$-terms by free variables. Clearly,
  this cannot increase the length $n$ of~\eqref{eq:p:h1}.
\end{proof}

\Section{Lower Bounds on Herbrand Disjunctions}

As noted in the proof of Theorem~\ref{t:1e}, the bound on the length of the
Herbrand disjunction depends only on the critical count of the initial proof.
This is in contrast to the bound we would obtain by the more standard approach
of cut-elimination and the mid-sequent theorem which depends on the length and
cut complexity of the original proof (see, e.g.,
\cite{Buss:1998,TroelSchwi:2000}).  In the case of the $\meps$-calculus, the
result concerns the relation between the critical count of a proof of $\exists
x\, E(x)$ in $\PCeps$ and the length of a Herbrand disjunction $\bigvee
E(t_i)$.  In the case of the sequent calculus and cut-elimination the result
concerns the relation between the length and cut complexity of a proof with
cut, and the length of a cut-free proof, which in turn determines the length
of a Herbrand disjunction obtained via the mid-sequent theorem.  In both
cases, the relationship is hyperexponential.  Statman~\cite{Statman:1979} and
Orevkov \cite{Orevkov:1982} showed that this bound is not just an artefact of
the particular cut-elimination procedure considered, but that proofs with cut
essentially have hyper-exponential speedup over cut-free proofs.  The question
may then be raised whether the same holds true of the $\meps$-calculus, i.e.,
whether the bound on the length of Herbrand disjunctions obtained in the first
$\meps$-theorem is tight.  Although we do not have a result quite as optimal
as Orevkov's in this regard, it can be shown that every $\meps$-elimination
procedure that yields Herbrand disjunctions must by hyperexponential.

We sketch the proof of such a lower bound theorem for the length of Herbrand
disjunctions. Recall that the \emph{Herbrand complexity}~$\HC(E)$ of a purely
existential formula $E \equiv \exists x_1 \ldots \exists x_n E'(x_1, \ldots,
x_n)$ is the length of the shortest $\lor$-expansion of~$E'(x_1, \ldots,
x_n)$.

\begin{theorem}
There is a sequence of formulas $E_k$ so that
\begin{enumerate}
\item for each $k$, there is a $\PCeps$-proof $\pi_k$ of $E_k$ with
  $\cc{\pi_k} \le c\cdot k$ (for some constant~$c$), but
\item $\HC(E_k) \ge 2_k^1$.
\end{enumerate}
\end{theorem}

We follow the presentation of Orevkov's Theorem in
\cite[\S6.11]{TroelSchwi:2000}. (Statman's result requires equality, but
Orevkov's does not.) Consider a language including a unary constant~$0$, a
unary function symbol~$S$ and a ternary relation~$R$, whose meaning is fixed
by the following axioms:
\begin{align*}
  \Hyp_1 & \equiv \forall x\, R(x,0,S(x)) \tkom \\
  \Hyp_2 & \equiv \forall y \forall x\forall z \forall z_1 (R(y,x,z) \land
  R(z,x,z_1) \impl R(y,S(x),z_1)) \tpkt
\end{align*}
Further, we define
\begin{equation*}
  C_k \equiv 
  \exists z_k \dots \exists z_0 (R(0,0,z_k) \land R(0,z_k,z_{k-1}) \land \dots
  \land R(0,z_1,z_0)) \tpkt
\end{equation*}
$R(n,m,k)$ expresses that $n + 2^m = k$, and $C_k$ expresses that $2_k^1$ is
defined. $E_k$ is the (purely existential) prefix form of $\Hyp_1 \land \Hyp_2
\impl C_k$.

\begin{lemma} 
  For every $k$, $\PCeps \vdash_{\pi_k} E_k$, where $\cc{\pi_k} = c\cdot k$
  (for some constant~$c$).
\end{lemma}

\begin{proof}
  $E_k$ is provable in the sequent calculus (alternatively, in natural
  deduction) using proofs (with cut) of length linear in~$k$,
  see~\cite{TroelSchwi:2000}. Proofs in the sequent calculus and in natural
  deduction can be translated into proofs in $\PCeps$ with linear increase in
  proof length. Moreover, as in the embedding lemma, only weak quantifier
  inferences ($\exists$I, $\forall$E) increase the critical count of the
  $\PCeps$-proof.  (We omit the details, which are routine.)
\end{proof}

This establishes part (1) of the theorem.  Orevkov's result concerns proof
lengths; we have to adapt the proof to Herbrand complexity.  We give a direct
proof of a lower bound on $\HC(E_k)$; the result can also be obtained using
techniques from proof complexity as in \cite{BaazLeitsch:1994:FI}. In order to
simplify the presentation, we will consider Herbrand \emph{sequents} of
$\Hyp_1, \Hyp_2 \seq C_k$ instead of Herbrand disjunctions, i.e., valid
sequents of the form $\Gamma_1, \Gamma_2 \seq \Delta$ where each formula in
$\Gamma_1$ is a substitution instance of $R(x, 0, S(x))$, each formula in
$\Gamma_2$ is a substitution instance of $R(y,x,z) \land R(z,x,z_1) \impl
R(y,S(x),z_1)$, and each formula in $\Delta$ is one of $R(0,0,z_k) \land
R(0,z_k,z_{k-1}) \land \dots \land R(0,z_1,z_0)$.  Then obviously
$\max\{\card{\Gamma_1}, \card{\Gamma_2}, \card{\Delta}\} \le \HC(E_k)$.  In
the following, $\bar n$ abbreviates $S^n(0)$, and $\Psi = \{\Hyp_1, \Hyp_2\}$.

\begin{lemma}\label{l:compute}
  Suppose $\Psi \seq R(\bar n, \bar m, \bar l)$ is valid. Then $l = n + 2^m$
  and for each Herbrand sequent $T \equiv (\Gamma_1, \Gamma_2 \seq R(\bar n,
  \bar m, \bar l))$ of $\Psi \seq R(\bar n, \bar m, \bar l)$, we have
  \begin{equation*}
    \{ R(\bar i, 0, S(\bar i)) : n \le i < n+2^m \} \subseteq \Gamma_1 \tpkt
  \end{equation*}
  In particular, $\card{\Gamma_1} \geq 2^m$.
\end{lemma}

\begin{proof}
  If $m=0$, then clearly the only possibility is $l = n+1$. Then any Herbrand
  sequent of $\Psi \seq R(\bar n, 0, S(\bar n))$ can be written as $R(\bar
  n,0,S(\bar n)), \Gamma' \seq R(\bar n, 0, S(\bar n))$ and satisfies the
  conditions.
  
  Suppose the result is established for $m$, and consider the case for $m+1$.
  Let $\Gamma_1,\Gamma_2 \seq R(\bar n, S(\bar m), \bar l)$ be any 
  Herbrand sequent of $\Psi \seq R(\bar n, S(\bar m), \bar l)$ with
  $\Gamma_1$ the instances corresponding to $\Hyp_1$ and $\Gamma_2$ those
  corresponding to~$\Hyp_2$.  As $R(\bar n, S(\bar m), \bar l)$ cannot follow
  from instances of $R(x, 0, S(x))$ alone, $\Gamma_2$ is nonempty and must
  contain a formula of the form
  \begin{equation*}
    R(\bar n,\bar m,\bar k) \land R(\bar k,\bar m,\bar l) \impl R(\bar n,
    S(\bar m), \bar l) \tkom
  \end{equation*}
  such that $T_1 \equiv (\Gamma_1, \Gamma_2 \seq R(\bar n,\bar m,\bar k))$ and
  $T_2 \equiv (\Gamma_1, \Gamma_2 \seq R(\bar k,\bar m,\bar l))$ are both
  valid.  That means that $T_1$ is a Herbrand sequent of $\Psi \seq R(\bar
  n,\bar m,\bar k)$ and $T_2$ one of $\Psi \seq R(\bar k,\bar m,\bar l)$.  The
  induction hypothesis applies, and it follows that $k= n+2^m$ and $l =
  k+2^m$, thus $l= n + 2^{m+1}$.  Further, $\Gamma_1$ must contain
  \begin{equation*}
    \{ R(\bar i, 0, S(\bar i)) : n \le i < n+2^m \} \cup
    \{ R(\bar i, 0, S(\bar i)) : n+2^m \le i < n+2^{m+1} \} \tkom
  \end{equation*}
  and the lemma follows.
\end{proof}

\begin{lemma}\label{l:minimal}
  Let $S$ be a sequent of the form $\Psi \seq \exists \bar z\, A(\bar z)$.
  Then every minimal Herbrand sequent of $S$ is of the form $\Gamma \seq
  \Delta$ with $\card{\Delta} = 1$.
\end{lemma}

\begin{proof}
  We use of the terminology and results of Chapter~XI of
  \cite{EbbinghausFlumThomas:1994}.  Suppose that $T \equiv \Gamma \seq
  \Delta$ is a Herbrand sequent of $S$ with $\Delta = A(\bar t_1), \ldots,
  A(\bar t_n)$.  Each formula in $\Gamma$ is Horn.  Thus the term model
  $\mathfrak{I}^\Gamma$ is a free model of $\Gamma$ (Corollary 2.5
  of~\cite{EbbinghausFlumThomas:1994}).  Since $\Gamma \seq \Delta$ is valid
  and $\mathfrak{I}^\Gamma \models \Gamma$, $\mathfrak{I}^\Gamma \models
  A(\bar t_1) \lor \ldots \lor A(\bar t_n)$.  Then there is an $i$ so that
  $\mathfrak{I}^\Gamma \models A(\bar t_i)$.  But $\mathfrak{I}^\Gamma$ is
  free. Hence, every model of $\Gamma$ is also a model of $A(t_i)$ and $\Gamma
  \seq A(t_i)$ is a Herbrand sequent of $S$.
\end{proof}

\begin{lemma}
  If $T\equiv(\Gamma_1, \Gamma_2 \seq \Delta)$ is a minimal Herbrand sequent
  of $\Psi \seq C_k$, then $\card{\Gamma_1} \geq 2_k^1$.
\end{lemma}

\begin{proof}
  By  Lemma~\ref{l:minimal}, $T$ is of the form
  \begin{equation} \label{HS}
    \Gamma_, \Gamma_2 \seq R(0,0, \bar n_k) \land R(0, \bar n_k, \bar n_{k-1}) \land \dots \land
    R(0, \bar n_1, \bar n_0) \tpkt
  \end{equation}
  (Note that by substituting $0$ for free variables, terms in a Herbrand
  sequent may always be brought into that form).  As~\eqref{HS} is valid, each
  of the sequents $\Gamma_1, \Gamma_2 \seq R(0,0, \bar n_k)$, \dots,
  $\Gamma_1,\Gamma_2 \seq R(0, \bar n_1, \bar n_0)$ is valid as well. Applying
  Lemma~\ref{l:compute} $(k-1)$-times, we see that $n_1 = 2_{k-1}^1$.  Since
  $\Gamma_1, \Gamma_2 \seq R(0, \bar n_1, \bar n_0)$ is a Herbrand sequent of
  $\Psi \seq R(0, \bar n_1, \bar n_0)$, $\Gamma_1$ contains the instances of
  $\Hyp_1$ given in Lemma~\ref{l:compute}, and $n_0 = 2^{n_1} = 2_k^1$. Hence,
  $\card{\Gamma} \ge 2_k^1$.
\end{proof}

\Section{The Second Epsilon Theorem}

The Second Epsilon Theorem is a generalization of the first.  It states that a
formula without $\meps$-terms provable in the extended predicate calculus is
already provable in the predicate calculus (without $\meps$-terms).  Its proof
utilizes no additional methods specific to the $\meps$-calculus beyond those
of the first $\meps$-theorem. 

\begin{definition}
  Suppose $A = {\sf Q}_1 x_1 \ldots {\sf Q}_n x_n\, B(x_1, \ldots, x_n)$ is a
  prenex formula.  Let $z_1$, \ldots, $z_l$ be all the $\forall$-quantified
  variables among $x_1$, \ldots, $x_n$, let $y_1$, \ldots, $y_m$ be all
  the $\exists$-quantified ones, and let $f_1$, \ldots, $f_l$ be new function
  symbols. The Herbrand normal form $A^H$ of $A$ is
  \begin{equation*}
    \exists y_1 \ldots \exists y_m\, 
    C(y_1, \ldots, y_m, t_1, \ldots, t_l) \tkom
  \end{equation*}
  where $t_j = f_j(y_1, \ldots, y_m)$.  
\end{definition}

\begin{lemma} \label{l:l5}
  Suppose $\PCeps \vdash A$. Then $\PCeps \vdash A^H$.
\end{lemma}

\begin{proof}
  Standard.
\end{proof}

\begin{theorem}[Second Epsilon Theorem]
  If $A$ is a formula of $L(\PC)$ and $\PCeps \vdash A$, then $\PC \vdash A$.
\end{theorem}

\begin{proof}
  For simplicity, we give the proof for prenex formula $A$ with a simple
  quantifier structure. The general results follows similarly. Assume $A$ has
  the form $\exists x \forall y \exists z\, B(x,y,z)$ with $B(x,y,z)$
  quantifier-free and only the indicated variables occur in $A$.  We apply
  Lemma~\ref{l:l5} to obtain a proof of $\exists x\exists z\, B(x, f(x), z)$.
  The extended first $\meps$-theorem now yields that there are $\meps$-free
  terms $r_{i},s_i$ so that
  \begin{equation} \label{eq:herbrand}
    \EC \vdash \bigvee_i B(r_{i}, f(r_i), s_i) \tpkt
  \end{equation}
  The idea is now to replace the $f(r_i)$ by new free variables $a_i$ and
  obtain from~\eqref{eq:herbrand}, that
  \begin{equation*}
    \bigvee_i B(r'_{i}, a_i, s'_i) 
  \end{equation*}
  is deducible in $\EC$.  Then the original prenex formula $A$ can be obtained
  if we employ suitable quantifier-shifting rules (deducible in $\PC$).
  
  Let $f(r_1), \ldots, f(r_p)$ denote terms occurring in the
  disjunction~\eqref{eq:herbrand}.  Let $p_i$ be the number of occurrences of
  $f$ in $f(r_{i})$.  We may assume that the disjunction is arranged so that
  that the sequence $f(r_1)$, \dots, $f(r_p)$ is ordered such that $p_i \le
  p_{i+1}$.  Now let $a_1$, \dots, $a_p$ be new free variables.  Replace each
  occurrence of $f(r_i)$ which does not occur as a subterm of another
  $f(r_{j})$ by $a_{i}$.  Then~\eqref{eq:herbrand} becomes
  \begin{equation}\label{eq:herbrand2}
    \bigvee_i B(r'_{i}, a_{i}, s'_i) \tkom
  \end{equation}
  Observe that $r'_i$ does not contain $a_j$ for $j \ge i$.  For if $r'_i$
  did contain $a_j$, then $r_i$ must contain $f(r_j)$, and $p_j < p_i$.  But we
  assumed that the disjunctions were ordered so that the sequence of $p_i$ was
  non-decreasing.
  
  It is easy to see that~\eqref{eq:herbrand2} is also a tautology, since pairs
  of equal atomic formulas remain pairs of equal atomic formulas.  Thus
  from~\eqref{eq:herbrand2} we obtain that
\begin{equation}\label{eq:herbrand3}
 \bigvee_i \exists z B(r'_{i}, a_i, z) \tkom
\end{equation}
by existentially generalizing on the terms $s'_i$.  
Now consider the last disjunct in ~\eqref{eq:herbrand3}. 
By the preceding observation, $a_p$ does not occur in any other disjunct,
or in $r_p'$.  Hence, in \PC, we may deduce from
\begin{equation*}
  \bigvee_{i=1}^{p-1} \exists z\, B(r'_{i}, a_{i}, z) \lor 
  \exists z\,B(r'_p, a_p, z) \tkom
\end{equation*}
the formula
\begin{equation*}
\bigvee_{i=1}^{p-1} \exists z\, B(r'_{i}, a_{i}, z) \lor 
     \forall y\exists z\,B(r'_p, y, z) \tpkt
\end{equation*}
Iterating these steps we eventually obtain a proof
of $A$ in $\PC$.
\end{proof}

\Section{Conclusion and Further Work}

The above proofs of the first and second $\meps$-theorem were formulated for
theorems of $\PCeps$. However, as indicated in the introduction, the theorems
remain valid in the presence of open (quantifier- and $\meps$-free) theories.

\begin{corollary}
Let $\Ax$ be a set of formulas without bound variables.
\begin{enumerate}
\item Let $E$ be a formula without bound variables (no quantifiers, no
  epsilons). If $\Ax \vdash E$ in $\PCeps$, then $\Ax \vdash E$ in $\EC$,
\item Let $\exists \overline{x} E(\overline{x})$ be a purely existential
  formula. If $\Ax \vdash \exists \overline{x} E(\overline{x})$ in $\PCeps$,
  then $\Ax \vdash \bigvee_i E(t_{i1}, \ldots, t_{in})$ for some $\meps$-free
  terms $t_{ij}$ in $\EC$
\item If $E$ is an $\meps$-free formula and $\Ax \vdash E$ in $\PCeps$, then
  $\Ax \vdash E$ in $\PC$.
\end{enumerate}
\end{corollary}

The discussion of the $\meps$-calculus given here is only a first step toward
a more comprehensive investigation of Hilbert's $\meps$-calculus.  The
gap in the upper and lower bound for the Herbrand complexity of theorems of
$\ECeps$ suggests that a stricter analysis or a refinement of the elimination
procedure for the first $\meps$-theorem is possible.  A more interesting and
pressing question, however, is the analysis of $\meps$-elimination procedures
for the $\meps$-calculus with equality.  The addition of equality to the
$\meps$-calculus is not as straightforward as it is in the predicate calculus,
and the first $\meps$-theorem is significantly more complicated if equality is
present than when it is not.  It remains a topic for future work.

In the introduction we claimed that encoding of quantifiers on the term level
using the $\meps$-operator may allow for a more condensed representation of
proofs.  Let us briefly sketch the reason for this. Since modus ponens is the
only inference rule in $\ECeps$, a formula $A^\meps$ is provable in $\ECeps$
iff there is a tautology of the form
\begin{equation}\label{eq:tautology}
  \bigwedge_{i,j} (B_i(t_j) \impl B_i(\eps{x}{B_i(x)})) \impl A^\meps \tkom
\end{equation}
Thus it suffices to find the critical formulas $B_i(t_j) \impl
B_i(\eps{x}{B_i(x)})$, i.e., the substitutions involved in the proof of~$A$,
such that~\eqref{eq:tautology} is a tautology.  This suggests that the
formalization of proofs is simpler in the $\meps$-calculus or at least that
proofs in the $\meps$-calculus can be represented more succinctly.

As pointed out above, the bound on the length of the Herbrand disjunctions
obtained using the first $\meps$-theorem depends only on the critical count of
the initial proof.  (If equality is present, however, the maximal rank of
critical formulas will also play a role.)  In other systems, such as the
sequent calculus, the number of critical formulas corresponds to the number of
weak quantifier inferences.  Standard methods for obtaining bounds on Herbrand
disjunctions ordinarily do not yield a bound only in the number of weak
quantifier inferences.  Consequently, the result obtained above is of
independent interest.  The change of input parameters is especially
significant when considered in conjunction with the above remarks on
formalizability. Standard methods usually only yield specific information
about Herbrand disjunctions---such as their length---if a complete formal
proof is available.  Within the $\meps$-calculus, we may weaken this
assumption, as provability witnessed by a tautology of the form
of~\eqref{eq:tautology} suffices.  This fact was successfully employed by
Kreisel in his ``unwinding'' of the proof of Littlewood's
theorem~\cite{Kreisel52-a}.

\bibliographystyle{studlog}
\begin{smallSL}
\newcommand{\bibfont}{\small} \setlength{\bibsep}{0pt}

\end{smallSL}

\AuthorAdressEmail{Georg Moser}{Computational Logic Group\\
Institute of Computer Science\\
University of Innsbruck\\
A--6020 Innsbruck, Austria}{georg.moser@uibk.ac.at}

\AdditionalAuthorAddressEmail{Richard Zach}{Department of Philosophy\\
University of Calgary\\
Calgary, AB T2N 1N4, Canada}{rzach@ucalgary.ca}

\label{l}

\end{document}